\documentclass[]{interact}

\usepackage{epstopdf}
\usepackage[caption=false]{subfig}
\usepackage{graphicx}
\usepackage{amsfonts}
\usepackage{amssymb}
\usepackage{amsmath}
\usepackage{enumerate} 
\usepackage{latexsym}
\usepackage[numbers,sort&compress]{natbib}
\bibpunct[, ]{[}{]}{,}{n}{,}{,}
\makeatletter
\def\NAT@def@citea{\def\@citea{\NAT@separator}}
\makeatother

\theoremstyle{plain}
\newtheorem{theorem}{Theorem}[section]
\newtheorem{lemma}[theorem]{Lemma}
\newtheorem{corollary}[theorem]{Corollary}
\newtheorem{proposition}[theorem]{Proposition}

\theoremstyle{definition}

\theoremstyle{remark}
\newtheorem{remark}{Remark}

\begin{document}


\title{Closed Range Composition Operators  on Hardy Spaces}

\author{
\name{Petros Galanopoulos\textsuperscript{a}\thanks{Contact Petros Galanopoulos Email: petrosgala@math.auth.gr} and Kostas Panteris\textsuperscript{b}\thanks{Contact Kostas Panteris  Email: mathp289@math.uoc.gr, kpanteris@yahoo.gr}}
\affil{\textsuperscript{a}Department of Mathematics, Aristotle University of Thessaloniki, 54124 Thessaloniki, Greece; 
\textsuperscript{b}Department of Mathematics and Applied Mathematics, University of Crete, University Campus Voutes, 70013 Heraklion, Greece }
}

\maketitle

\begin{abstract}
If $\varphi$ is an analytic self-map of the open unit disc $\mathbb{D}$ in the complex plane, the composition operator $C_{\varphi}$ on spaces of analytic functions is defined as $C_{\varphi}(f) = f \circ \varphi$. In this paper we prove two equivalent conditions for the composition operator $C_{\varphi}$ on $H^{p}$, $0 < p < \infty$, to have closed range. One condition is a nontrivial extension, with a proof following a different approach, of a condition proved by Zorboska (1994) for $H^2$. The other condition is a trivial extension of a condition proved by Cima, Thomson, and Wogen (1974) again for $H^2$ and we include it in this work for the sake of completeness. We also prove, as an application of our previous results and the theory of Aleksandrov-Clark measures, that if $\varphi$ is an inner function then $C_{\varphi}$ has closed range on $H^{p}$, 
\end{abstract}


\section{Introduction and Preliminaries}
 Let $\mathbb{D}$ denote the open unit disk in the complex plane, $\mathbb{T}$ the unit circle, $A$ the normalized area Lebesgue measure in $\mathbb{D}$ and $m$ the normalized length Lebesgue measure in $\mathbb{T}$. For $0 < p < \infty$ the Hardy space $H^{p}$ is defined as the set of all analytic functions in $\mathbb{D}$ for which
\[
\sup \limits_{0\leq r<1} \int_0^{2\pi} \!\!  \vert f(r e^{i\theta}) \vert^{p} dm(\theta) < +\infty
\] 
and the corresponding norm in $H^{p}$ is defined by
\[
\Vert f  \Vert_{H^{p}}^{p}  = \sup \limits_{0 \leq r<1} \int_0^{2\pi} \!\!  \vert f(r e^{i\theta}) \vert^{p} dm(\theta).
\]

In this work we will mainly make use of the following equivalent norm (see Hardy-Stein identities in \cite{Pavlovic}, pages 58-59):
\begin{equation}
\Vert f  \Vert_{H^{p}}^{p}  = \vert f(0)  \vert^{p} + \iint \limits_{\mathbb{D}}  \vert f(z)\vert^{p-2} \vert f^{\prime}(z) \vert^{2} \log\frac{1}{\vert  z \vert} dA(z).
\end{equation}

If $\varphi$ is a non-constant analytic self-map of the unit disk $\mathbb{D}$, then the composition operator $C_{\varphi}: H^{p} \rightarrow H^{p}$ is defined as $C_{\varphi}(f) = f \circ \varphi$ and the Nevanlinna counting function $N_{\varphi}$ is defined as  
\begin{equation}\label{Nev_eq}
    N_{\varphi}(w) =
    \begin{cases}
      \sum \limits_{\varphi(z) = w} \log\frac{1}{\vert  z \vert}, & \text{if}\ w \in \varphi(\mathbb{D}) \setminus \lbrace \varphi(0)  \rbrace \\
      0, & \text{otherwise}\ 
    \end{cases}
\end{equation}

Let $\rho(z,w)$ denote the pseudo-hyberbolic distance between $z,w \in \mathbb{D}$,
\[
\rho(z,w) = \Big\vert \frac{z - w}{1 - \overline{z}w} \Big\vert,
\] 
and $D_{\eta}(a)$ denote the pseudo-hyberbolic disk of center $a \in \mathbb{D}$ and radius $\eta<1$:
\[
D_{\eta}(a) = \lbrace z \in \mathbb{D}: \rho(a,z) < \eta  \rbrace.
\]

In the following, $C$ denotes a positive and finite constant which may change from one occurrence to another. Moreover, by writing
$K(z) \asymp L(z)$ for the non-negative quantities $K(z)$ and $L(z)$ we mean that $K(z)$ is comparable to $L(z)$ if $z$ belongs to a specific set: there are positive constants 
 $C_{1}$ and $C_{2}$ independent of $z$ such that
\[
C_{1} K(z) \leq L(z) \leq C_{2} K(z).
\]
 
 In \cite{CTW} and \cite{Z} the case of closed range composition operators in Hardy space $H^{2}$ is studied.

 In \cite{CTW}, Cima, Thomson, and Wogen  gave an equivalent condition for $C_{\varphi}:H^{2} \rightarrow H^{2}$ to have closed range
that depends only on the behavior of the function $\varphi$ on the boundary $\mathbb{T}$ of the open unit disk $\mathbb{D}$.
 First, they extend $\varphi$ in $\mathbb{T}$ as $ \varphi(\zeta)= \lim \limits_{r \rightarrow 1} \varphi(r \zeta)$, since it is well known that this limit exists for almost every $\zeta\in\mathbb{T}$ with respect to Lebesgue measure $m$. Then they define the measure $\nu_{\varphi}$ on  Borel sets
$E \subset \mathbb{T}$ by  
\begin{equation}\label{ni_measure_phi}
\nu_{\varphi}(E)= m(\varphi^{-1}(E)).
\end{equation}
The measure $\nu_{\varphi}$ is absolutely continuous with respect to Lebesgue
measure $m$ on $\mathbb{T}$ and its Radon-Nikodym derivative $\frac{d\nu_{\varphi}}{dm}$
is in $L^{\infty}(\mathbb{T})$. 
\begin{theorem}[Cima, Thomson, and Wogen]\label{Cima}
$C_{\varphi}:H^{2} \rightarrow H^{2}$ has closed range if and only if the Radon-Nikodym derivative $\frac{d\nu_{\varphi}}{dm}$ is essentially bounded away from zero.
\end{theorem}
 
 In \cite{Z}, Zorboska proved a criterion for $C_{\varphi}$ to have closed range on $H^{2}$ based upon properties of $\varphi$ on
 pseudo-hyberbolic disks. She defines the function 
\[
\tau_{\varphi}(z) = \frac{N_{\varphi}(z)}{\log\frac{1}{\vert  z \vert}}
\]
for $z \in \varphi(\mathbb{D}) \setminus \varphi(0)$ and, for $c>0$, the set 
\begin{equation}\label{G_set}
G_{c} = \lbrace z \in  \mathbb{D}:  \tau_{\varphi}(z)  > c \rbrace.
\end{equation} 
 
\begin{theorem}[Zorboska]\label{Zorb}
$C_{\varphi}:H^{2} \rightarrow H^{2}$ has closed range if and
only if there exist constants $c>0$, $\delta > 0$ and $\eta \in (0,1)$ such that the set $G_{c}$ 
satisfies
\[
A(G_{c} \cap D_{\eta}(a)) \geq \delta A(D_{\eta}(a))
\]
for all $a \in \mathbb{D}$.
\end{theorem}

\section{Main result}
We are going to prove that the results of theorems \ref{Cima} and \ref{Zorb} hold, not only for $H^{2}$, but for every $H^{p}$, $0 < p < \infty$. Here is the main result. 

\begin{theorem}\label{main_result}
Let $0<p<\infty$. Then the following are equivalent:
\begin{enumerate}
\item[(i)] $C_{\varphi}: H^{p} \rightarrow H^{p}$ has closed range.
\item[(ii)] The Radon-Nikodym derivative $\frac{d\nu_{\varphi}}{dm}$ is essentially bounded away from zero.
\item[(iii)] There exist $c>0$, $\delta > 0$ and $\eta \in (0,1)$ such that the set $G_{c}$ satisfies
\[
A(G_{c} \cap D_{\eta}(a)) \geq \delta A(D_{\eta}(a))
\]
for all $a \in \mathbb{D}$.
\end{enumerate}
\end{theorem}

\par The equivalence $(i) \Leftrightarrow (ii)$ is actually theorem \ref{Cima}. Our proof is the same and we include it for the sake of completeness. The equivalence $(i) \Leftrightarrow (iii)$ is theorem \ref{Zorb}. As it is not clear whether Zorboska's proof for $p=2$ works for every $p>0$, we prove this equivalence following a different approach. Namely, we use Hardy-Stein identities (see \cite{Pavlovic}, pages 58-59) for one of the directions, and reverse Carleson measures (see \cite{Nicolau}) and pull-back measures (see \cite{Pull_back}) for the converse.
\par We will make use of the following lemma \ref{pb2} and theorem \ref{pb1}, proved in \cite{Pull_back}, as well as theorem \ref{rev}. The case $p>1$ of theorem \ref{rev} is proved in \cite{Nicolau}.

Let $\Delta = \frac{4\partial^{2}}{\partial z \partial \overline{z}}$ be the usual Laplacian and, for $\zeta \in \mathbb{T}$ and $0 < h < 1$, let $W(\zeta, h)$ be the usual Carleson square
\[
W(\zeta, h) = \lbrace  z \in \overline{\mathbb{D}}: 1-h < \vert z \vert \leq 1, \vert \arg(z \overline{\zeta})  \vert \leq \pi h \rbrace.
\]

We will also make use of the measure $m_{\varphi}$ defined on Borel sets
$E \subset \overline{\mathbb{D}}$ by 
\begin{equation}\label{measure_phi}
m_{\varphi}(E)= m(\varphi^{-1}(E)\cap\mathbb{T}).
\end{equation}
Actually, $\nu_{\varphi}$ defined in \eqref{ni_measure_phi} is the restriction of $m_{\varphi}$ on $\mathbb{T}$. 

\begin{lemma}\label{pb2}
For every $g \in C^{2}(\mathbb{C})$ we have
\[
\iint \limits_{\overline{\mathbb{D}}} g(z) dm_{\varphi}(z) = g(\varphi(0)) + \frac{1}{2} \iint \limits_{\mathbb{D}} \Delta g(w) N_{\varphi}(w) dA(w).
\]
\end{lemma}

\begin{theorem}\label{pb1}
For $0 < c < \frac{1}{8}$, $\zeta \in \mathbb{T}$  and   $0 < h < (1 - \vert \varphi(0)   \vert) / 8$, we have 
\[
\sup \limits_{z \in W(\zeta, h) \cap \mathbb{D} } N_{\varphi}(z) \leq  \frac{100}{c^{2}} m_{\varphi}( W(\zeta, (1 + c)h) ).
\]
\end{theorem}

\begin{theorem}\label{rev}
Let $0 < p < \infty$ and let $\mu$ be a positive measure in $\overline{\mathbb{D}}$. Then the following assertions are equivalent.
\begin{enumerate}[(i)]
\item[(i)] There exists $C > 0$ such that for every $f \in H^{p} \cap C(\overline{\mathbb{D}})$,
\[
\iint \limits_{\overline{\mathbb{D}}} \vert f(z) \vert^{p} d\mu(z) \geq C \Vert f  \Vert_{H^{p}}^{p}.
\]
\item[(ii)] There exists $C > 0$ such that for every $\lambda \in \mathbb{D}$
\[
\iint \limits_{\overline{\mathbb{D}}} \vert K_{\lambda}(z) \vert^{p} d\mu(z) \geq C.
\]
where, for $p>1$, $k_{\lambda}(z) = \frac{1}{1 - \overline{\lambda} z}$ is the reproducing kernel in $H^{p}$ and $K_{\lambda} = \frac{k_{\lambda}}{\Vert  k_{\lambda} \Vert_{H^{p}}}$ is its normalised version, and, for $0 < p \leq 1$, we have $K_{\lambda}(z) = \frac{1 - \vert \lambda \vert^{2}}{(1 - \overline{\lambda} z )^{(p+1)/p}} $.
\item[(iii)] There exists $C > 0$ such that for $0<h<1$ and $\zeta\in\mathbb{T}$ we have
\[
\mu(W(\zeta,h)) \geq C h.
\]
\item[(iv)] There exists $C > 0$ such that the Radon-Nikodym derivative of $\mu\vert_{\mathbb{T}}$ (the restriction of $\mu$ on $\mathbb{T}$) with respect to $m$ is bounded below by $C$.
\end{enumerate}
\end{theorem}
 \begin{remark} 
 As we have already mentioned, the case $p>1$ of theorem \ref{rev} is proved in \cite{Nicolau}. For the case $0 < p \leq 1$, we have just to observe that, with the choice of $K_{\lambda}(z) = \frac{1 - \vert \lambda \vert^{2}}{(1 - \overline{\lambda} z )^{(p+1)/p}} $ in assertion $(ii)$ of theorem \ref{rev}, then the proof of all assertions of the theorem, as described  in \cite{Nicolau}, works also for the case $0 < p \leq 1$.
\end{remark}

In the proof of theorem \ref{main_result} we will also make use of a theorem of D. Luecking in \cite{Luecking81} which we restate here as
\begin{theorem}\label{Lue}
Let $\tau$ be a non-negative, measurable, bounded function in $\mathbb{D}$ and for $c>0$ let $G_{c} = \lbrace z \in  \mathbb{D}:  \tau(z)  > c \rbrace $. Then the following assertions are equivalent.
\begin{enumerate}
\item[(i)] There exists $C>0$ such that 
\[
\iint \limits_{\mathbb{D}} \vert f^{\prime}(z) \vert^{2}  \tau(z)  \log\frac{1}{\vert z \vert} dA(z) 
\geq C \iint \limits_{\mathbb{D}} \vert f^{\prime}(z) \vert^{2} \log\frac{1}{\vert z \vert} dA(z)
\]
for every $f \in H^{2}$.
\item[(ii)] There exist $C>0$ and $c>0$ such that 
\[
\iint \limits_{G_{c}} \vert f^{\prime}(z) \vert^{2} \log\frac{1}{\vert z \vert} dA(z) 
\geq C \iint \limits_{\mathbb{D}} \vert f^{\prime}(z) \vert^{2} \log\frac{1}{\vert z \vert} dA(z)
\]
for every $f \in H^{2}$.
\item[(iii)] There exist $c>0$, $\delta > 0$ and $\eta \in (0,1)$ such that
\[
A(G_{c} \cap D_{\eta}(a)) \geq \delta A(D_{\eta}(a))
\]
for all $a \in \mathbb{D}$.
\end{enumerate}
\end{theorem}

In the following we will also use of the following non-univalent change of variable formula (see \cite{Shapiro}, section 4.3). If $g$ is a measurable, non-negative function in $\mathbb{D}$, we have
\begin{equation}\label{change_v}
\iint \limits_{\mathbb{D}} g(\varphi(z)) \vert  \varphi^{\prime}(z) \vert^{2} \log \frac{1}{\vert z \vert} dA(z) = 2 \iint \limits_{\mathbb{D}}  g(w) N_{\varphi}(w) dA(w).
\end{equation}

\begin{proof}[Proof of theorem \ref{main_result}.]
Without loss of generality we may assume that $\varphi(0)=0$ and $f(0)=0$ for all functions $f$ in $H^p$.\\
$(i) \Rightarrow (iii)$ If $C_{\varphi}$ has closed range then there exist $C>0$ (we may suppose $C < 2$) such that for every $f \in H^{p}$ we have
\[
\Vert  C_{\varphi}f \Vert_{H^{p}}^{p}  \geq C \Vert  f  \Vert_{H^{p}}^{p} 
\]
i.e.
\begin{align*}
\iint \limits_{\mathbb{D}} \vert f(\varphi(z)) \vert^{p-2} \vert f^{\prime}(\varphi(z)) \vert^{2} \vert & \varphi^{\prime}(z) \vert^{2} \log\frac{1}{\vert  z \vert} dA(z) \\
 &\geq C \iint \limits_{\mathbb{D}} \vert f(z) \vert^{p-2} \vert f^{\prime}(z) \vert^{2} \log\frac{1}{\vert  z \vert} dA(z).
\end{align*}
By \eqref{change_v} we have
 \begin{align*}
 \iint \limits_{\mathbb{D}} \vert f(w) \vert^{p-2} \vert f^{\prime}(w) \vert^{2} & N_{\varphi}(w)  dA(w) \\
 & \geq C \iint \limits_{\mathbb{D}} \vert f(z) \vert^{p-2} \vert f^{\prime}(z) \vert^{2} \log\frac{1}{\vert  z \vert} dA(z)
\end{align*}
i.e.
\begin{align}\label{interm}
 \iint \limits_{\mathbb{D}} \vert f(w) \vert^{p-2} \vert f^{\prime}(w) \vert^{2} &\tau_{\varphi}(w) \log\frac{1}{\vert  w \vert} dA(w) \nonumber\\
 &\geq C \iint \limits_{\mathbb{D}} \vert f(z) \vert^{p-2} \vert f^{\prime}(z) \vert^{2} \log\frac{1}{\vert  z \vert} dA(z).
\end{align}
Let $f \in H^{p}$ with $f(z) \neq 0$ for every $z \in \mathbb{D}$. We define the analytic function
\[
g(z) = f(z)^{p/2}.
\]
Obviously, $g \in H^{2}$ and $g(z) \neq 0$ for every $z \in \mathbb{D}$. Then from \eqref{interm} we have
\[
\iint \limits_{\mathbb{D}} \vert g^{\prime}(w) \vert^{2} \tau_{\varphi}(w) \log\frac{1}{\vert  w \vert} dA(w)
 \geq C \iint \limits_{\mathbb{D}} \vert g^{\prime}(z) \vert^{2} \log\frac{1}{\vert  z \vert} dA(z)
\]
and, because of Luecking's theorem \ref{Lue},
\[
\iint \limits_{G_{c}} \vert g^{\prime}(w) \vert^{2} \log\frac{1}{\vert  w \vert} dA(w)
 \geq C \iint \limits_{\mathbb{D}} \vert g^{\prime}(z) \vert^{2} \log\frac{1}{\vert  z \vert} dA(z)
\]
or, equivalently,
\begin{equation}\label{interm2}
\iint \limits_{G_{c}} \vert g^{\prime}(w) \vert^{2} (1 - \vert  w \vert^{2}) dA(w)
 \geq C \iint \limits_{\mathbb{D}} \vert g^{\prime}(z) \vert^{2} (1 - \vert  z \vert^{2}) dA(z).
\end{equation}
for every $g \in H^{2}$ with $g(z) \neq 0$ for every $z \in \mathbb{D}$.\\ 
Let $a \in \mathbb{D}$. Choosing $g \in H^{2}$ such that  $\vert g^{\prime}(z)  \vert^{2} = \frac{(1 - \vert  a \vert^{2})^{3}}{\vert 1 - \overline{a}z  \vert^{6}}$ and with $g(z) \neq 0$ for every $z \in \mathbb{D}$, we get

\begin{align}\label{basic_rel}
\iint\limits_{G_{c} \cap D_{\eta}(a) } \frac{(1 - \vert  a \vert^{2})^{3}}{\vert 1 - \overline{a}z  \vert^{6}}  (1 - \vert  z \vert^{2}) dA(z) 
 \geq &\iint \limits_{G_{c} } \frac{(1 - \vert  a \vert^{2})^{3}}{\vert 1 - \overline{a}z  \vert^{6}} (1 - \vert  z \vert^{2}) dA(z)\nonumber \\
&- \iint \limits_{\mathbb{D} \setminus D_{\eta}(a) } \frac{(1 - \vert  a \vert^{2})^{3}}{\vert 1 - \overline{a}z  \vert^{6}} (1 - \vert  z \vert^{2}) dA(z).
\end{align}
From \eqref{interm2} we get
\begin{equation}\label{interm3}
\iint \limits_{G_{c}} \frac{(1 - \vert  a \vert^{2})^{3}}{\vert 1 - \overline{a}z  \vert^{6}} (1 - \vert  z \vert^{2}) dA(z) \geq C \iint \limits_{\mathbb{D}} \frac{(1 - \vert  a \vert^{2})^{3}}{\vert 1 - \overline{a}z  \vert^{6}} (1 - \vert  z \vert^{2}) dA(z)=\frac{C}{2}.
\end{equation}
We have
\[
\iint \limits_{\mathbb{D}} (1 - \vert  z \vert^{2}) dA(z) = \frac{1}{2}.
\]
and we choose $\eta \in (0, 1)$ such that
\[
\iint \limits_{\vert z \vert < \eta} (1 - \vert  z \vert^{2}) dA(z) \geq \frac{1}{2} - \frac{C}{4}= \frac{2-C}{4}>0.
\]
By the change of variable $z = \frac{w - a}{1 - \overline{a}w} = \psi_{a}(w)$ we get
\[
\iint \limits_{D_{\eta}(a)} (1 - \vert  \psi_{a}(w) \vert^{2}) \frac{(1 - \vert  a \vert^{2})^{2}}{\vert 1 - \overline{a}w  \vert^{4}} dA(w) \geq \frac{2-C}{4}
\]
and hence
\begin{equation}\label{interm4}
\iint \limits_{\mathbb{D} \setminus D_{\eta}(a)} (1 - \vert  \psi_{a}(w) \vert^{2}) \frac{(1 - \vert  a \vert^{2})^{2}}{\vert 1 - \overline{a}w  \vert^{4}} dA(w) \leq \frac{C}{4}.
\end{equation}
Combining \eqref{basic_rel}, \eqref{interm3} and \eqref{interm4}, we find
\[
\iint \limits_{G_{c} \cap D_{\eta}(a) } \frac{(1 - \vert  a \vert^{2})^{3}}{\vert 1 - \overline{a}z  \vert^{6}} (1 - \vert  z \vert^{2}) dA(z) \geq \frac{C}{4}.
\]
Using the fact that if $z \in D_{\eta}(a)$ then $(1 - \vert  a \vert^{2}) \asymp (1 - \vert  z \vert^{2}) \asymp \vert 1 - \overline{a}z \vert $, we get
\[
\frac{C^{\prime} A(G_{c} \cap D_{\eta}(a))}{(1 - \vert  a \vert^{2})^{2}} \geq \frac{C}{4}
\]
and, finally,
\[
A(G_{c} \cap D_{\eta}(a))\geq \delta A(D_{\eta}(a)),
\]
$(iii) \Rightarrow (i)$ We consider the measure $m_{\varphi}$ as defined in \eqref{measure_phi} and we will show that (iii) of theorem \ref{main_result} implies (iii) of theorem \ref{rev} with $\mu=m_{\varphi}$.\\ 
Now we consider $\zeta \in \mathbb{T}$ and $0<h<1$ and the corresponding Carleson square $W(\zeta, h)$. Having in mind to apply theorem \ref{pb1}, we take $c=\frac 1{16}$ and $h'=\frac h8$. Then there exists $a \in W(\zeta, h')$ so that
\[
D_{\eta}(a) \subset W(\zeta, h') \subset W(\zeta, (1+c)h') \subset W(\zeta, h),\qquad 1-\vert a \vert^2\geq Ch,
\] 
where $C$ depends upon $\eta$. We have $A(G_{c} \cap D_{\eta}(a)) > \delta A(D_{\eta}(a))$ and hence $G_{c} \cap D_{\eta}(a) \neq \emptyset$. Let $b \in G_{c} \cap D_{\eta}(a)$. Then $1 - \vert b \vert^{2} \geq Ch$ and $N_{\varphi}(b) \geq c \log\frac{1}{\vert b \vert}$ (since $b \in G_{c}$). Applying theorem \ref{pb1} (recalling that $\varphi(0) = 0$), we find

\begin{align*}
m_{\varphi}(W(\zeta, h)) & \geq  m_{\varphi}(W(\zeta, (1+c)h')) \geq C \sup \limits_{z \in W(\zeta, h)\cap\mathbb{D}} N_{\varphi}(z) \geq C N_{\varphi}(b)\\
& \geq C \log\frac{1}{\vert b \vert}   \geq C (1 - \vert b \vert^{2}) \geq Ch.
\end{align*}
Therefore we get (iii) of theorem \ref{rev} which is equivalent to (i) of the same theorem, with $\mu=m_{\varphi}$. Now we take any $f$ which is analytic in a disk larger than $\mathbb{D}$ and so that $f(0)=0$ and we use (i) of theorem \ref{rev} together with lemma \ref{pb2} to find
\begin{align*}
\Vert  C_{\varphi}f \Vert_{H^{p}}^{p} & = \iint \limits_{\mathbb{D}} \vert f(\varphi(z)) \vert^{p-2} \vert f^{\prime}(\varphi(z)) \vert^{2} \vert \varphi^{\prime}(z) \vert^{2} \log\frac{1}{\vert  z \vert} dA(z)\\ 
& = \iint \limits_{\mathbb{D}} \vert f(w) \vert^{p-2} \vert f^{\prime}(w) \vert^{2} N_{\varphi}(w)  dA(w)\\
& =  \frac{1}{p} \iint \limits_{\mathbb{D}} \Delta(\vert f \vert^{p})  N_{\varphi}(w)  dA(w) \\ 
& = C \iint \limits_{\overline{\mathbb{D}}} \vert f(w) \vert^{p} dm_{\varphi}(w)\\
& \geq C \Vert f  \Vert_{H^{p}}^{p},
\end{align*}
Now, if $f$ is the general function in $H^{p}$ with $f(0)=0$, we apply the result to the functions $f_{r}$, $0<r<1$, defined by $f_{r}(z) = f(r z)$,  $z \in \mathbb{D}$, and we take the limit as $r\rightarrow 1-$. Therefore $C_{\varphi}$ has closed range.\\
$(i) \Rightarrow (ii)$  Let's suppose that $C_{\varphi}$ has closed range on $H^{p}$ and $E \subset \mathbb{T}$. For $n \in \mathbb{N}$ we choose $f_{n} \in H^{p}$ (theorem 4.4, page 63, in \cite{Garnett}) such that
\begin{equation*}
    \vert f_{n}(\zeta) \vert^{p} =
    \begin{cases}
      1, & \text{if}\ \zeta \in E \\
      \frac{1}{2^{n}}, & \text{if}\ \zeta \in \mathbb{T}\setminus E 
    \end{cases}
\end{equation*}
Then $\Vert  C_{\varphi}f_{n} \Vert_{H^{p}}^{p} \geq C \Vert  f_{n} \Vert_{H^{p}}^{p}$ and hence 
\begin{align*}
 m(\varphi^{-1}(E)) + \frac{1}{2^{n}} m(\mathbb{T} \setminus \varphi^{-1}(E))
 \geq C m(E)  + C \frac{1}{2^{n}} m(\mathbb{T}\setminus E).
\end{align*}
Taking limit as $n \rightarrow +\infty$, we get $m(\varphi^{-1}(E)) \geq C m(E)$, i.e.
\[
\nu_{\varphi}(E) \geq C m(E).
\]
Thus the Radon-Nikodym derivative $\frac{d\nu_{\varphi}}{dm}$ is bounded below by $C$.\\
$(ii) \Rightarrow (i)$ Let's suppose that the Radon-Nikodym derivative $\frac{d\nu_{\varphi}}{m}$ is bounded below by $C.$ For $\lambda > 0$ we consider the set
\[
E_{f}(\lambda) = \big\lbrace e^{i\theta}: \vert f(e^{i\theta}) \vert > \lambda  \big\rbrace.
\]
Then, 
\[
m\big(E_{f\circ\varphi}(\lambda)\big) = \nu_{\varphi}\big( E_{f}(\lambda) \big) \geq C m\big( E_{f}(\lambda) \big)
\]
for all $\lambda>0$, and finally
\begin{align*}
&\Vert  C_{\varphi}f \Vert_{H^{p}}^{p}  =  \int \limits_{0}^{2\pi} \vert f(\varphi(e^{i\theta})) \vert^{p} dm(\theta) = \int \limits_{0}^{+\infty} p \lambda^{p-1} m\big( E_{f\circ\varphi}(\lambda) \big) d\lambda  \\
& \geq C \int \limits_{0}^{+\infty} p \lambda^{p-1} m\big( E_{f}(\lambda) \big) d\lambda = C \int \limits_{0}^{2\pi} \vert f(e^{i\theta}) \vert^{p} dm(\theta) = C \Vert  f \Vert_{H^{p}}^{p}.
\end{align*} 
Hence $C_{\varphi}$ has closed range.
\end{proof}
Now we are going to show that a result regarding Besov type spaces, due to M. Tjani \cite{Tjani}, can be extended by using theorem \ref{main_result}. If $p>1, \alpha>-1$, the Besov type space $B_{p, \alpha}$ is the space of analytic functions $f$ on $\mathbb{D}$ such that
\begin{equation*}
\Vert f \Vert_{B_{p, \alpha}}^{p} = |f(0)|^p + \iint \limits_{\mathbb{D}}  \vert f^{\prime}(z) \vert^{p} (1 - \vert  z \vert^{2})^{\alpha} dA(z) < +\infty
\end{equation*}
If $\varepsilon>0$, the set $G_{\varepsilon, p, \alpha}$ is defined as
\begin{equation*}
G_{\varepsilon, p, \alpha} = \Big\lbrace w \in \mathbb{D} : \frac{N_{p, \alpha}(w, \varphi)}{(1 - \vert w \vert^{2})^{\alpha}}>\varepsilon \Big\rbrace
\end{equation*}
where
\begin{equation*}
N_{p, \alpha}(w, \varphi) = \sum \limits_{\varphi(z) = w} \vert \varphi^{\prime}(z)  \vert^{p-2} (1 - \vert z \vert^{2})^{\alpha}
\end{equation*}
In \cite{Tjani} the following results are proved.
\begin{theorem}\label{Tjani_theo}
(Theorem 5.2 in \cite{Tjani}). For $p>2$, the operator $C_{\varphi}: B_{p, p-1} \rightarrow  B_{p, p-1}$ has closed range if and only if there exists an $\varepsilon>0$ so that $G_{\varepsilon, p, p-1}$ satisfies a reverse Carleson condition, which means that there exist $\varepsilon>0$, $\delta>0$ and $\eta \in (0,1)$ such that $A(G_{\varepsilon, p, p-1} \cap D_{\eta}(z)) > \delta A(D_{\eta}(z))$ for every $z \in \mathbb{D}$.
\end{theorem}
\begin{corollary}\label{Tjani_corol}
(Corollary 5.3 in \cite{Tjani}). Let $p>2$. If $C_{\varphi}$ has closed range on $B_{p, p-1}$ then $C_{\varphi}$ has closed range on $H^{2}$.
\end{corollary}

By using theorem \ref{main_result} we can extend corollary \ref{Tjani_corol} and get the following result.
\begin{corollary}\label{extended_corol}
Let $p>2$. If $C_{\varphi}$ has closed range on $B_{p, p-1}$ then $C_{\varphi}$ has closed range on every $H^{q}$, $0< q < \infty$.
\end{corollary}
\begin{proof}
If $C_{\varphi}$ has closed range on $B_{p, p-1}$ then, by theorem \ref{Tjani_theo}, we have that there exist $\varepsilon>0$, $\delta>0$ and $\eta \in (0,1)$ such that $A(G_{\varepsilon, p, p-1} \cap D_{\eta}(z)) > \delta A(D_{\eta}(z))$ for every $z \in \mathbb{D}$. Moreover, we have for the set $G_{\varepsilon}$, as defined in \eqref{G_set}, that $G_{\varepsilon} = G_{\varepsilon,2,1} $ and also $G_{\varepsilon, p, p-1} \subset G_{\varepsilon,2,1}$. Finally we have that 
\begin{equation*}
A(G_{\varepsilon} \cap D_{\eta}(z)) > A(G_{\varepsilon, p, p-1} \cap D_{\eta}(z)) > \delta A(D_{\eta}(z))
\end{equation*}
so, by theorem \ref{main_result}, $C_{\varphi}$ is closed range on every $H^{q}$, $0< q<\infty$.
\end{proof}

\section{Regarding inner functions}

It is well-known, as a theorem of Herglotz, that if $u$ is a non-negative harmonic function in $\mathbb{D}$, then there exists a unique positive Borel measure $\mu$ such that 
\[
u(z) = \int \limits_{\mathbb{T}} \frac{1 - \vert z  \vert^{2}}{\vert \zeta - z  \vert^{2}} d\mu(\zeta).
\]
Now, if $\varphi:\mathbb{D} \rightarrow \mathbb{D}$ is an analytic function and $\alpha\in\mathbb{T}$, then the function 
\[
\Re\Big( \frac{\alpha + \varphi(z)}{\alpha - \varphi(z)}  \Big) = \frac{1 - \vert \varphi(z)  \vert^{2}}{\vert \alpha - \varphi(z)  \vert^{2}}
\]
is positive and harmonic in $\mathbb{D}$ and from Herglotz's theorem we have
\[
\frac{1 - \vert \varphi(z)  \vert^{2}}{\vert \alpha - \varphi(z)  \vert^{2}} = \int \limits_{\mathbb{T}} \frac{1 - \vert z  \vert^{2}}{\vert \zeta - z  \vert^{2}} d\mu_{\alpha}(\zeta)
\]
for a unique Borel measure $\mu_{\alpha}$ in $\mathbb{T}$. The measures $\mu_{\alpha}$, $\alpha\in\mathbb{T}$, are called Alexandrov-Clark measures. Now we have the Lebesgue decomposition 
\[
d\mu_{\alpha} = h_{\alpha} dm + d\mu_{\alpha}^{s}, \hspace{2mm} h_{\alpha}\in L^{1}(\mathbb{T}),\hspace{2mm} \mu_{\alpha}^{s} \perp m. 
\]
It is well known (see \cite{Cima_Alexandrov}, pages 204-208) that the total variation of $\mu_{\alpha}$ is given by 
\begin{equation}\label{total_variation}
\Vert \mu_{\alpha} \Vert = \frac{1 - \vert \varphi(0)  \vert^{2}}{\vert \alpha - \varphi(0)  \vert^{2}},
\end{equation}
that the absolutely continuous part $h_{\alpha} dm$ is carried by the set $\lbrace \zeta\in\mathbb{T}:\vert \varphi(\zeta)\vert<1\rbrace$,
\begin{equation}\label{absolutely_continuous}
h_{\alpha}(\zeta) = \frac{1 - \vert \varphi(\zeta)  \vert^{2}}{\vert \alpha - \varphi(\zeta)  \vert^{2}}
\end{equation}
and that the singular part $d\mu_{\alpha}^{s}$ is carried by the set $\lbrace \zeta\in\mathbb{T}: \varphi(\zeta) = \alpha \rbrace$.
From \eqref{total_variation} and \eqref{absolutely_continuous}, we see that
\begin{align}\label{singular part}
\Vert \mu_{\alpha}^{s} \Vert & = \Vert \mu_{\alpha} \Vert - \int \limits_{\mathbb{T}} h_{\alpha}(\zeta)dm(\zeta)\nonumber\\ 
&=\frac{1 - \vert \varphi(0) \vert^{2}}{\vert \alpha - \varphi(0)  \vert^{2}}- \int \limits_{\mathbb{T}} \frac{1 - \vert \varphi(\zeta)  \vert^{2}}{\vert \alpha - \varphi(\zeta)  \vert^{2}} dm(\zeta).
\end{align}

The result which connects the Alexandrov-Clark measures with $C_{\varphi}$ having closed range is a proposition due to K. Luery in \cite{Luery} (page 56).
\begin{proposition}[Luery]\label{prop_Luery}
For $m$-a.e $\alpha\in\mathbb{T}$,
\begin{equation}\label{basic_Luery_rel}
\frac{d\nu_{\varphi}}{dm}(\alpha) = \Vert \mu_{\alpha}^{s} \Vert. 
\end{equation}
\end{proposition}
Now we can prove the following result.
\begin{proposition}\label{inner_hardy_prop}
Let $0<p<\infty$. If $\varphi:\mathbb{D} \rightarrow \mathbb{D}$ is inner then $C_{\varphi}:H^{p} \rightarrow H^{p}$ has closed range.
\end{proposition}
\begin{proof}
From proposition \ref{prop_Luery}, from relation \eqref{singular part}, and from the fact that $\vert \varphi(\zeta)  \vert=1$ for $m$-a.e $\zeta\in\mathbb{T}$ when $\varphi$ is inner, we get 
\begin{align}
\frac{d\nu_{\varphi}}{dm}(\alpha) & = \frac{1 - \vert \varphi(0) \vert^{2}}{\vert \alpha - \varphi(0)  \vert^{2}} - \int \limits_{\mathbb{T}} \frac{1 - \vert \varphi(\zeta)  \vert^{2}}{\vert \alpha - \varphi(\zeta)  \vert^{2}} dm(\zeta)\nonumber\\
& = \frac{1 - \vert \varphi(0) \vert^{2}}{\vert \alpha - \varphi(0)  \vert^{2}}\geq\frac{1 - \vert \varphi(0) \vert^{2}}{4}>0
\end{align}
for $m$-a.e $\alpha\in\mathbb{T}$.

From part (ii) of theorem \ref{main_result} we conclude that $C_{\varphi}$ has closed range.
\end{proof}
\section*{Acknowledgements}
Many thanks to Prof. Michael Papadimitrakis for discussions about the mathematical content of this paper. His contribution was essential in order for it to take its final form.

%
%



%
%

\end{document}